# Role of the clay lenses within sandy aquifers in the migration pathway of infiltrating DNAPL plume: A numerical investigation

Shikhar Nilabh, Fidel Grandia

Amphos 21 Consulting S.L., c/Veneçuela 103, 08019 Barcelona, Spain

## Abstract

The use of numerical based multi-phase fluid flow simulation can significantly aid in the development of an effective remediation strategy for groundwater systems contaminated with Dense Non Aqueous Phase Liquid (DNAPL). Incorporating the lithological heterogeneities of the aquifer into the model domain is a crucial aspect in the development of robust numerical simulators. Previous research studies have attempted to incorporate lithological heterogeneities into the domain; however, most of these numerical simulators are based on Finite Volume Method (FVM) and Finite Difference Method (FDM) which have limited applicability in the field-scale aquifers. Finite Element Method (FEM) can be highly useful in developing the field-scale simulation of DNAPL infiltration due to its consistent accuracy on irregular study domain, and the availability of higher orders of basis functions.

In this research work, FEM based model has been developed to simulate the DNAPL infiltration in a hypothetical field-scale aquifer. The model results demonstrate the effect of meso-scale heterogeneities, specifically clay lenses, on the migration and accumulation of Dense Non Aqueous Phase Liquid (DNAPL) within the aquifer. Furthermore, this research provides valuable insights for the development of an appropriate remediation strategy for a general contaminated aquifer.

## 1 Introduction

The increase in industrial activity over the past decades has led to a significant rise in anthropogenic contamination of groundwater, with Dense Non Aqueous Phase Liquids (DNAPLs) being a common form of these contaminants in aquifers. Due to its higher density than water, DNAPLs tend to migrate to the deeper levels of the groundwater system (Li *et al.*, 2022; Luciano *et al.*, 2010; Pankow and Cherry, 1996). These contaminants have a high persistence in the groundwater system, with residence times ranging from a few years to decades (Heron *et al.*, 2016; Kavanaugh *et al.*, 2003). Despite the development of various groundwater remediation techniques, their efficiency is often hindered by a poor understanding of the spatio-temporal evolution of DNAPLs in the subsurface (Kavanaugh *et al.*, 2003; Mohan Kumar and Mathew, 2005, Nilabh et al., 2022). The sparsity of data in field studies, coupled with the heterogeneity of hydrogeological properties, often results in an incomplete information of the DNAPL distribution within the aquifer. In such scenarios, modeling techniques can be utilized to estimate the source zone architecture of DNAPL by numerically reconstructing its migration pathway (Kamon *et al.*, 2004; McLaren *et al.*, 2012). Previous modeling studies in the literature have consistently shown the strong dependence of DNAPL migration pathway on heterogeneities (Gerhard and Kueper, 2003; He *et al.*, 2022; Jakobs *et al.*, 2003; Zheng *et al.*, 2015). Meso-scale heterogeneities, such as clay lenses, play a significant role in controlling the accumulation or diversion of DNAPL flow based on its hydrogeological properties and morphology (Ayral-Çınar and Demond, 2020; Reynolds and Kueper, 2001). The formation of a secondary source of contamination within an aquifer can occur either on top of low permeability layers or as residual saturation along migration pathways, as illustrated in Figure 1. As a result, a comprehensive understanding of the behavior of Dense Non Aqueous Phase Liquids (DNAPL) in the presence of clay lenses is crucial for the development of an effective groundwater remediation plan.

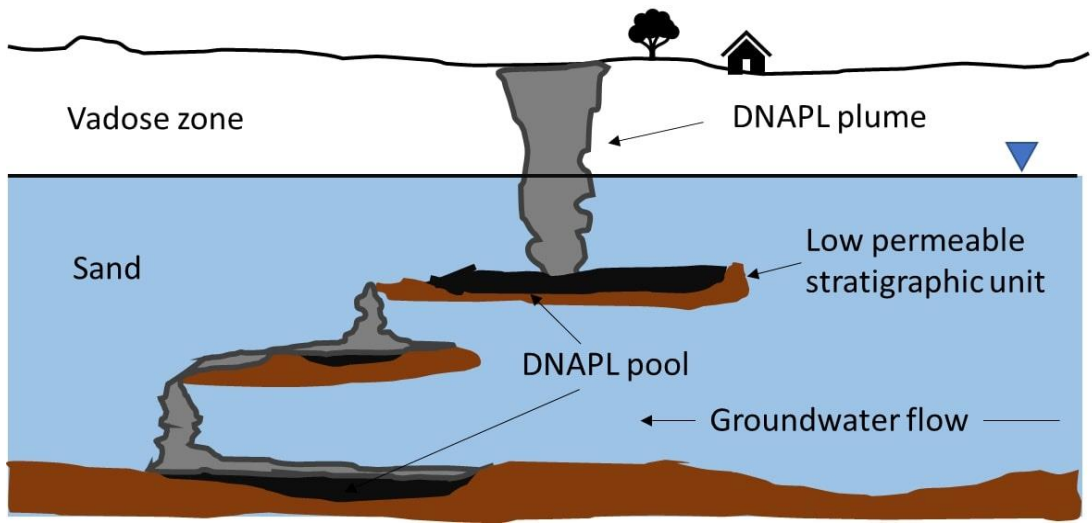

*Figure 1*: Scheme of DNAPL infiltration in a typical groundwater system, resulting in the formation of a complex DNAPL source zone architecture

Previous research studies have developed a variety of modeling tools for simulating the migration pathway of Dense Non Aqueous Phase Liquids (DNAPL). These tools have been implemented in codes such as SECUReFoam (Icardi *et al.*, 2022), DuMu$^x$ (Koch *et al.*, 2021), Tough2 v2 (Pruess *et al.*, 2011), UTCHEM (CPGE, 2000), STOMP (Nichols *et al.*, 1997) and VENT2D (Benson, 1994). A common challenge in the development of numerical formulations for two-phase flow in heterogeneous media is addressing discontinuities in hydrogeological properties. A mathematical form of interface condition has been proposed in the literature to explicitly define the discontinuities in the fluid´s pressure and saturation profile of water at the lithological interface (Helmig et al., 2006; Helmig, 1997). In recent years, several models have included this interface condition in their numerical formulations for two-phase flow (Kou and Sun, 2010; Papafotiou *et al.*, 2010; Tran *et al.*, 2020). While there are several models based on finite difference (Blagojevic *et al.*, 2022; Oostrom *et al.*, 2007; Suk *et al.*, 2019) or finite volume methods (Alali *et al.*, 2021; Kou and Sun, 2010; Tran *et al.*, 2020) for DNAPL infiltration in an aquifer, the flexibility of geometry and intrinsic boundary adaptation of these models is limited. A model based on Finite Element Method (FEM) offers more flexibility with the numerical environment and thus, can be more robust and favorable for field-scale aquifers (Li *et al.*, 2020, Yu and Li, 2019). However, FEM-based numerical models for simulating DNAPL infiltration in a field-scale scenario are still uncommon in the published literature (Yu and Li, 2019).

In this article, a FEM based formulation is applied to simulate DNAPL infiltration in a hypothetical field-scale aquifer representing a general polluted site with Perchloroethylene (PCE). The hypothetical aquifer is characterised by meso-scale heterogeneities represented by several low permeable clay units. Specifically, this modeling exercise is intended to demonstrate the application of developed model on deriving relevant predictive insight such as 1) the migratory pathway for DNAPL phase, 2) the time period of infiltration 3) the DNAPL accumulation and branch formation 4) effect of clay lenses on the migration pathway, 5) the DNAPL source zone architecture, 6) the amount of DNAPL pool and ganglia formed, 7) the numerical effect of the implementation of the Interface Condition and capillary forces, 8) the mass conservation of the DNAPL in the model domain, 9) the insights for real field management.

## 2 Methodology

### 2.1 Mathematical formulation

The mathematical formation governing the two-phase flow in the groundwater system has been briefly reviewed. The common equation for fluid flow in continuum scale is Darcy's law, which defines the relationship between the hydraulic gradient and the fluid flow rate and is given by Equation 1:

$$\boldsymbol{Q = -KA(\nabla H)} \qquad (eq.1)$$

Where Q [$L^3 \cdot T^{-1}$] is the volumetric fluid flow rate through a cross section of area A [$L^2$], K is the hydraulic conductivity [$L \cdot T^{-1}$] represented by Equation 2:

$$\boldsymbol{K = \frac{k\rho g}{\mu}} \qquad (eq.\ 2)$$

$k$ [L$^2$] is permeability of the porous media. $\rho$ [M·L$^{-3}$] is density of the fluid, $\mu$ [M·L$^{-1}$T$^1$] is the viscosity of the same fluid and g [L·T$^{-2}$] is the gravity constant,

H [L] is the hydraulic head defined as Equation 3:

$$H = (\frac{\nabla p}{\rho g} + z) \quad (eq.\ 3)$$

Where $p$ [M·L$^{-1}$T$^2$] is the pressure of the fluid, z [L] is the elevation accounting for the potential energy due to gravity in the z direction.

The Darcy velocity is defined as the flow rate of the fluid per unit cross section area and given by Equations 4 and 5:

$$u = \frac{Q}{A} = -K(\nabla H) \quad (eq.\ 4)$$

$$u = -\frac{k}{\mu}(\nabla p - \rho g) \quad (eq.\ 5)$$

Where $u$ is the Darcy velocity of the fluid.

The above equations assume that the flow is laminar within an incompressible porous medium (Grant, 2005). A general multiphase flow is represented by their mass conservation equation as shown in Equation 6 (Abriola and Pinder, 1985; Corapcioglu and Baehr, 1987):

$$\frac{\partial(\theta \rho_\alpha S_\alpha)}{\partial t} = -\nabla \cdot (\rho_\alpha u_\alpha) + \rho_\alpha q_\alpha \qquad \alpha = w, o \quad (eq.\ 6)$$

Where $\theta$ is the porous media, $S_\alpha$ is the saturation of fluid phase $\alpha$, $q_\alpha$ is the volumetric source or sink term, and $k_{r\alpha}$ is the relative permeability of phase $\alpha$ for the same fluid $\alpha$.

Therefore for the water-DNAPL system, the set of equations becomes (eq. 7, 8 and 9):

$$\frac{\partial(\theta S_w)}{\partial t} = \nabla \cdot \left(k \frac{k_{rw}}{\mu_w}(\nabla p_w - \rho_w g)\right) + q_w \quad (eq.\ 7)$$

$$\frac{\partial(\theta S_n)}{\partial t} = \nabla \cdot \left(k \frac{k_{rn}}{\mu_n}(\nabla p_n - \rho_n g)\right) + q_n \quad (eq.\ 8)$$

$$S_n + S_w = 1 \quad (eq.\ 9)$$

Where the subscript $w$ and $n$ with the variables are water and DNAPL, respectively.

Mathematically, there are four dependent variables: pressure of water, pressure of DNAPL, saturation of water and saturation of DNAPL. Since there are only 3 equations, one more equation is needed. An empirical relationship between capillary pressure and saturation of water is used for this purpose which empirically relates these two variables (eq. 10)

$$p_c(S_w) = p_n - p_w \qquad (eq.\ 10)$$

### 2.1.1 Two-phase flow in heterogeneous media: Interface condition

The Phase Pressure Saturation with Interface Condition (PPSIC) has been postulated in previous studies to simulate the flow dynamics at the interface of different sand layers (Helmig *et al.*, 2002; Van Duijn *et al.*, 1995). The theory is based on a capillary pressure equilibrium condition at the interface of two different sands for the DNAPL plume to infiltrate into a heterogeneous porous media. For an aquifer free of DNAPL, the capillary pressure obtained using the Brooks-Corey relationship is the same as the entry pressure (Brooks and Corey, 1964). The different entry pressures for different sands cause discontinuity in the capillary pressure. For the flow of DNAPL from higher permeability to lower permeability sands, the DNAPL encounters a relatively higher entry pressure at the plume front.

The interface condition for DNAPL infiltration through different sandy layers incorporates the capillary pressure as a quasi-primary variable in the numerical formulation. For the DNAPL plume infiltration through two sand layers (Sand A to Sand B) with different hydrogeological properties can be expressed as Equations 11, 12, 13 and 14:

For non-equilibrium capillary condition ($p_c$ at Sand A $<$ $p_c$ at Sand B)

$$S_{e_1} \geq S_e^* \qquad (eq.\ 11)$$

$$S_{e_2} = 1 \qquad (eq.\ 12)$$

For Capillary Equilibrium condition ($p_c$ at Sand A $=$ $p_c$ at Sand B)

$$S_{e_1} < S_e^* \qquad (eq.\ 13)$$

$$p_{d_1} S_{e_1}^{\frac{-1}{\lambda_1}} = p_{d_2} S_{e_2}^{\frac{-1}{\lambda_2}} \qquad (eq.\ 14)$$

Where $p_{d_1}$ and $p_{d_2}$ are the entry pressure of sand A and sand B, respectively. $S_{e_1}$ and $S_{e_2}$ are the effective saturation of wetting phase at the interface for sand A and sand B, respectively. $\lambda_1$ and $\lambda_2$ is the pore distribution coefficient for sand A and sand B, respectively.

$S_e^*$ is the threshold value of effective water saturation below which the capillary equilibrium is maintained. This occurs when the accumulation of the DNAPL on the top of clay lens replaces water from the pores below the threshold value to attain the capillary equilibrium condition. In the case of equilibrium condition, the different entry pressure and pore distribution coefficient, following Equation 14 yields different values of effective water saturation at the interface. This process is reflected as a discontinuity in the DNAPL saturation curve. Until the capillary pressure equilibrium is not attained, the movement of DNAPL is stopped. As the capillary pressure is attained once the accumulated DNAPL is sufficient enough to overcome threshold pressure, the movement resumes. At this time, the capillary pressure equilibrium is maintained but the saturation of DNAPL remains the same. Developing two-phase formulation without considering interface can fail to correctly simulate the flow dynamics in an aquifer with different sand layers. Therefore, for developing the field scale fluid flow simulation, the interface condition is included in the numerical formulation.

## 2.2 Comsol software package

The numerical model in this work is built on the Comsol Multiphysics software package (comsol 2021), which uses Galerkin based Finite Element method for solving Partial Differential Equations (PDEs). The software´s capability of coupling diverse physical phenomena as well as customizing the mesh and solver for each problem make it suitable for developing the robust two-phase flow model. The software uses MUltifrontal Massively Parallel sparse (MUMPS) Solver, a direct method for solving the linear system of equation obtained from Finite element analysis. The model represents the flow dynamics at a macroscopic scale which allows a continuum approach of defining parameter and variables in Comsol. The numerical formulation used for modeling has been already verified with the benchmark models in the literature (Nilabh, 2021)

## 2.3 Numerical model implementation

### 2.3.1 Model domain and boundaries

A hypothetical aquifer of dimension 50 m ×15 m × 15m is considered for simulating the DNAPL infiltration (Figure 2). The simulated domain represents a typical sandy aquifer with meso-scale heterogeneities in the form of low permeable clay lenses. The infiltration area of DNAPL has been assigned on the top of the model domain as shown by the orange arrow.

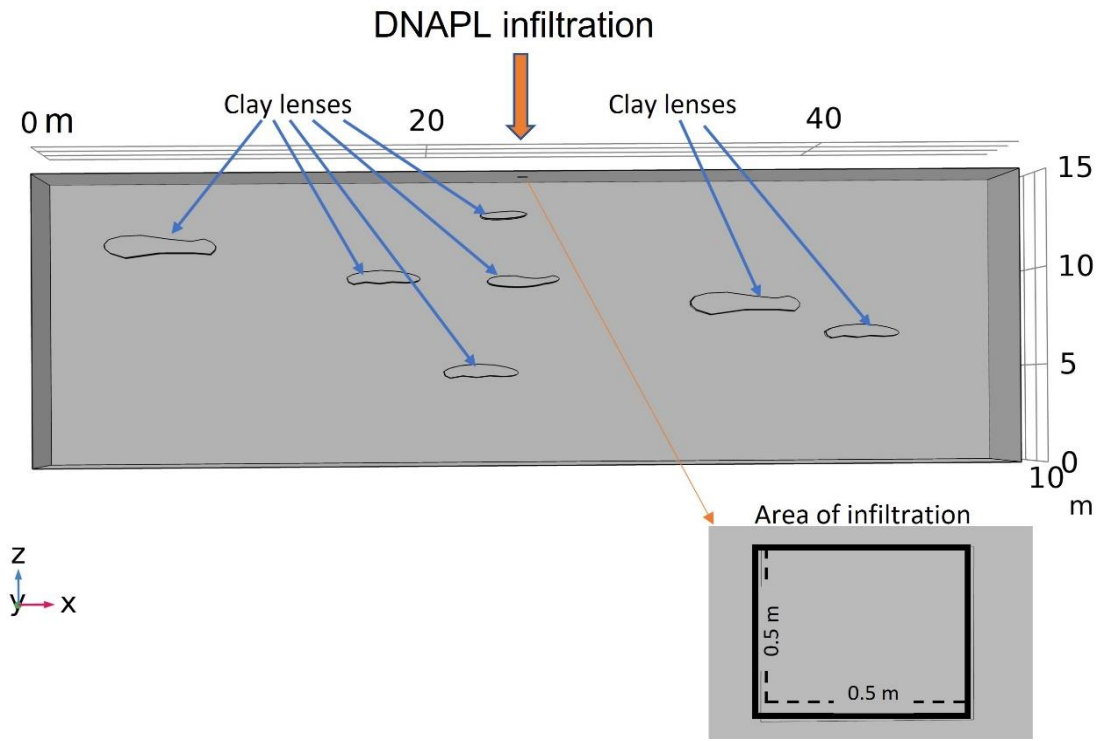

*Figure 2*: *Model domain used for simulation of two-phase flow. The domain simulates homogeneous sand with randomly placed clay lenses acting as meso-scale heterogeneities. An area with the dimension of 0.5m X 0.5 m is considered at the top of model domain through which the DNAPL infiltration occurs.*

One of the key questions concerning the migration of DNAPL in sandy aquifers is the impact of clay-bearing lenses or clay plugs, which are frequent in many examples of sandy layers in literature (Ayral-Çınar and Demond, 2020; Guo *et al.*, 2021; Hastings and Lee, 2021). To see such an impact, six clay lenses have been included in the model (Figure 2). These clay lenses have a random spatial location, geometry, with lengths varying from 1 to 2 m. These clay lenses are hydraulically of much lower permeability compared to the surrounding sand ($1.5 \times 10^{-10}$ m$^2$ vs. $5 \times 10^{-14}$ m$^2$), and, thus, the DNAPL flow is expected to be hindered at the sand-clay interface. The lithological properties used for the hypothetical aquifer are taken partially from a real field case of Innovation Garage, which is a research site located in an industrial area near the Copenhagen city centre (Denmark). This site is being used by The Capital Region of Denmark to study the migration of organochloride solvents in soils and sandy aquifers, and to test efficient remediation methodologies (The Capital Region Denmark, 2017). Due to unavailability of the complete dataset from the field case, the missing dataset has been taken from the literature (Kueper *et al.*, 1991). This approach ensures that the simulation result lies well within the realm of possible DNAPL infiltration event in a general field setting. Finally, the model assumes that the water table, recharge rate and all the intrinsic properties does not change through the simulation period.

Table 1: Hydrogeological properties used in the two-phase model of DNPL infiltration within the hypothetical aquifer.

| Property | Sand | Clay |
|---|---|---|
| Permeability | $1.5\times10^{-10}$ m$^2$ | $5\times10^{-14}$ m$^2$ |
| Porosity | 0.3 | 0.2 |
| Water residual saturation | 0.098 | 0.19 |
| Oil residual saturation | 0.01 | 0.008 |
| Entry pressure | 1323 Pa | 4500 Pa |
| Soil distribution index | 3.86 | 3.51 |

The decoupled water pressure ($P_w$) and DNAPL Saturation ($S_n$) formulation, also known as Implicit Pressure, Explicit Saturation (IMPES) formulation has been selected as a numerical scheme (Kou and Sun, 2010). The primary variable are the pressure of water and the saturation of the DNAPL, and the numerical framework requires initial and boundary conditions for these primary variables. The water in the aquifer has been considered to be hydrostatic and free of DNAPL at t = 0. The density of simulated DNAPL has been considered to be 1630 kg·m$^3$ to represent perchloroethylene (PCE), one of the organochlorides that infiltrated in Innovation Garage site. The hypothetical aquifer is considered to have a slow dynamics of groundwater flow (<10cm/day) and thus have negligible effect on the DNAPL migration pathway. This is in accordance with the previous studies, which reported negligible groundwater flow on DNAPL migration, given the natural flow velocity is less than 10 cm/day (Zheng *et al.*, 2015).

As the aquifer is assumed to be free of DNAPL at the onset of simulation, the initial DNAPL saturation should be zero. However, a minimal DNAPL saturation of 0.001 has been considered in the initial state of aquifer to make the achievement of the numerical convergence criteria. An arbitrarily defined infiltration rate of 0.375g/s of DNAPL is conceptualised to enter from the area of the infiltration placed at top of the model domain (Figure 2). A total of 500 Kg of DNAPL entry is simulated during the first 15 days of infiltration

### 2.3.2 Meshing

The mesh is designed considering the expected flow behavior of DNAPL in the simulated domain. A heterogeneous mesh geometry is used for the spatial discretization, which consists of a total of 2,304,156 tetrahedral and hexahedral elements. Figure 3a shows the geometry and size of the mesh elements used for the spatial discretization of the domain. In this figure, a symmetric half of the domain is removed for the better visualization of the mesh elements. The mesh size and geometry has been adjusted such that the spatial gradient of the state variables can be reasonably captured in the discretization process (Figure 3). Here, the mesh size of a finite element signifies the length of the longest edge of that element. As the infiltration of DNAPL is implemented from the central top, a relatively smaller mesh size is used throughout the central vertical section of the model. For the region away from the central part of the domain, a relatively coarser element size has been used. The maximum element size used is 10.2 meters in length, while minimum size used is 1.7 centimeters After several iteration of modelling the DNAPL infiltration, an optimized mesh configuration is used at the expected DNAPL flow path. The Figure 3b shows the mesh elements used along the expected flow path of the simulated DNAPL plume for the upper half of model domain. Finer elements in the range of 2 centimeters to 15 centimeters are used at the periphery of the clay lens to capture the flow behavior of DNAPL at the lithological interfaces. For the lower half of the model domain, a hexahedral system of mesh elements is used along the expected flow path of the DNAPL (Figure 3c). The finer mesh is generated at the base of the bedrock to capture the accumulation of DNAPL plume at the bottom of the domain. For the time discretization, an adaptive time discretization method is used. The initial time step was considered to be $1\times10^{-7}$ seconds, but according to the convergence criteria, the time discretization becomes coarser for an efficient computational time. The simulation of the DNAPL infiltration in a field scale domain has been run for 1500 hours.

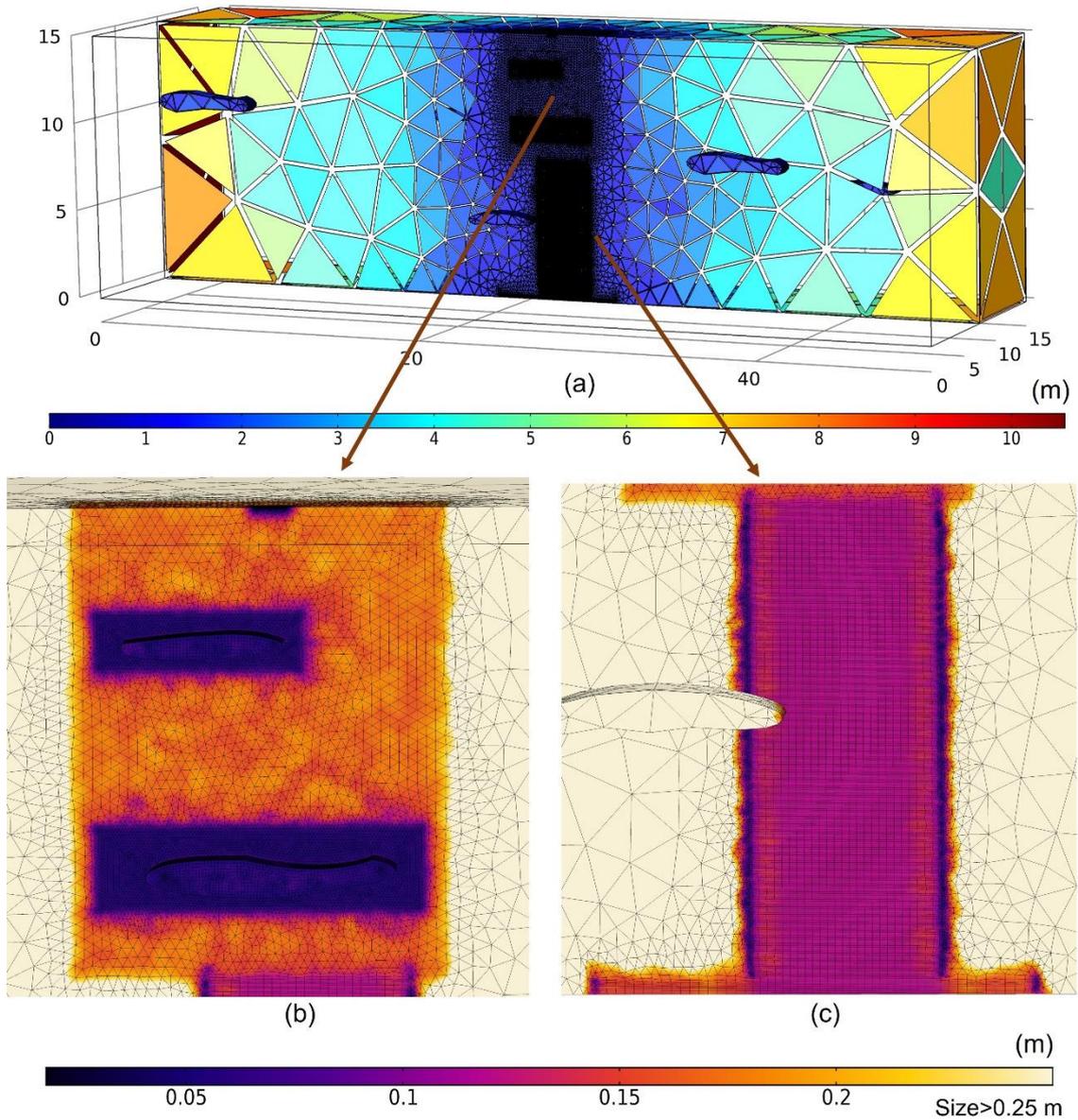

*Figure 3*: *(a) Spatial discretization of the study domain for finite element approximation with tetrahedral and hexahedral elements of size ranging from 0.017 to 10.2 m. (b,c) Enlarged image of mesh elements used along the expected DNAPL flow path*

## 3   Results

For understanding the results from the model, the DNAPL infiltration has been studied along two cross section planes XY and YZ, as shown in Figure 4. These planes crosscut the aquifer symmetrically into two halves to capture the DNAPL movement dynamics.

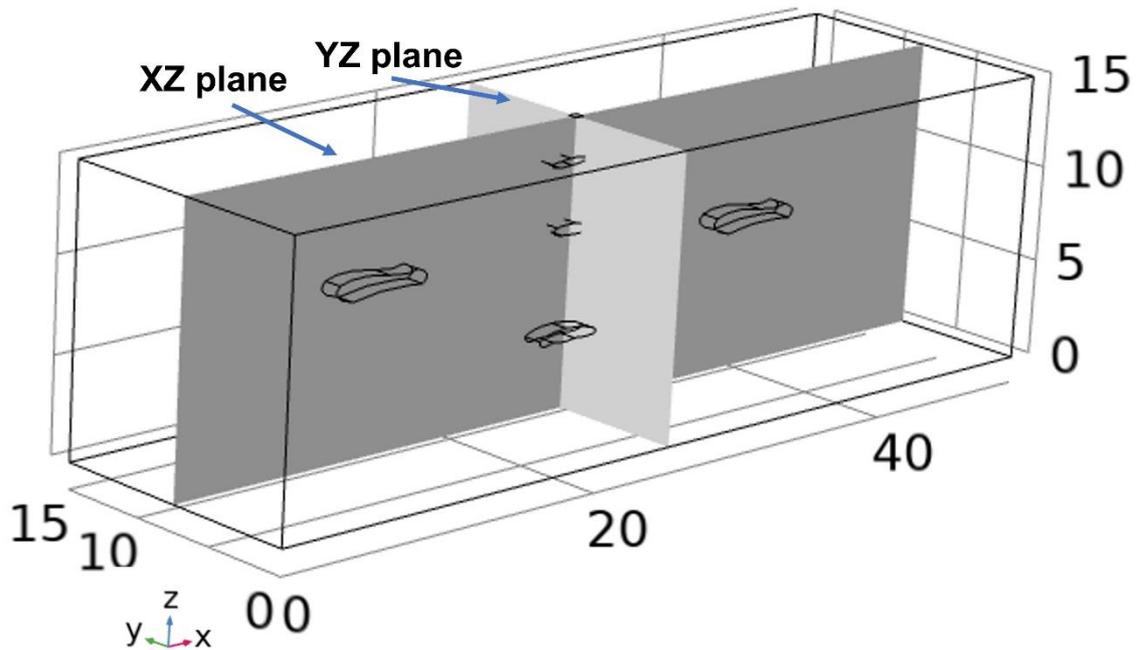

*Figure 4*: Two reference cross section planes aligned in XZ and YZ directions and dividing the simulated aquifer into symmetrical halves. These planes are used for result visualization.

## 3.1 Prediction of DNAPL infiltration

The model predicts that DNAPL infiltration takes 42 days to reach the bottom of the aquifer. The migration path is quite unique and dependent on both aquifer characteristics as well as DNAPL properties. The model also predicts that DNAPL enters the aquifer system with a saturation of 0.15 with a vertical infiltration of 1.5 meter after the end of the first day (Figure 5). The horizontal extent of DNAPL predicted by model is 0.5 m which is nearly the same as the dimension of inlet. Therefore, the predominant migration pathway aligns with the vertical direction. Thus, the result implies that the gravitational forces are the dominant for the DNAPL migration compared to viscous and capillary forces. This illustrates that for the field with homogeneous sand, the DNAPL source architecture is likely to have a simple vertical geometry with minimal lateral migration of the plume.

The predicted DNAPL plume infiltrates with steady state saturation for 34 hours until it meets the clay lens occurring at a depth of 1.9 meters. At the clay lens, the model predicts accumulation of DNAPL (Figure 5c-d), leading to the increase of DNAPL saturation from 0.15 to 0.41 in the following 18 hours. This accumulation causes lateral migration of DNAPL both along the width and the length of the clay lens. The model predicts that the DNAPL lateral movement along the width is much slower compared to that along the length of DNAPL. This highlights the formation of preferred channels along which most of the DNAPL moves along the length of the clay. Such preferred channels forms because the lateral DNAPL plume movement is estimated to reach the end of the clay lens comparatively faster than the plume along the width. This leads to re-occurrence of vertical movement of simulated DNAPL plume in the sandy aquifer, while its movement is still lateral along the width. The formation of preferred channels of DNAPL plume leads to the formation of unique DNAPL plume geometry. After 3 days, the model predicts the DNAPL plume infiltration up to 3.5 meters deep.

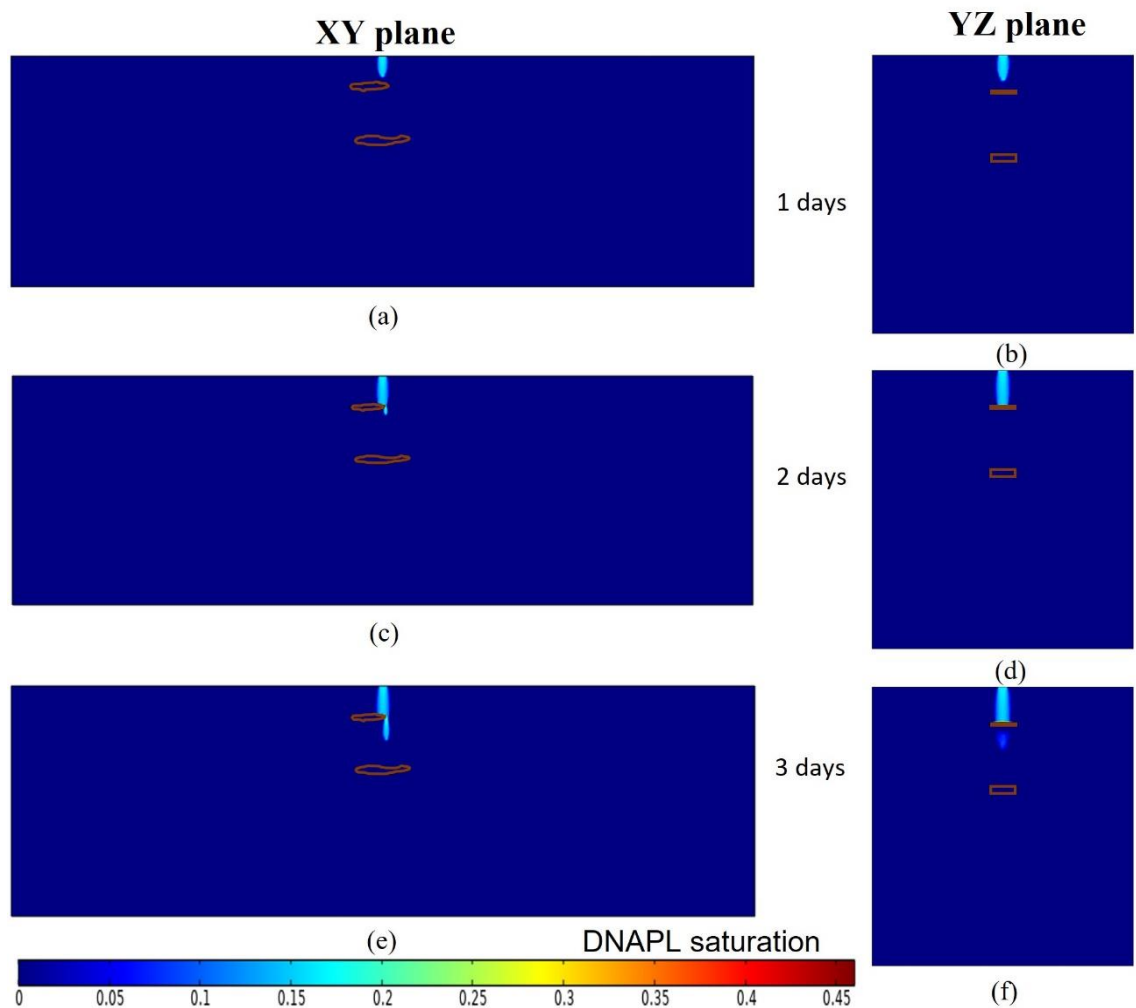

*Figure 5*: *Model prediction for DNAPL saturation in the study domain along the XZ plane (a) and YZ plane (b) after 1 day. After 2 days the model predicts further infiltration along the XZ plane (c) and YZ plane (d). After 3 days, when the source condition is removed, the model estimation shows the DNAPL profile along the XZ plane (e) and YZ plane (f). The clay lenses are outlined with brown color for a better visualization.*

The model also predicts that the multiple branches of DNAPL would form around the clay in the aquifer. After 5.5 days of infiltration, the DNAPL plume is predicted to settle down on the lower clay lens occurring at 9.5 meters deep. The Figure 6 shows the model prediction for DNAPL migration for 6 days, 10 days and 15 days, respectively. The accumulation of DNAPL over the clay lens is predicted to result at DNAPL saturation of 0.45, until the plume forms the branch. The model also predicts that the DNAPL plume branching along the width of the lower clay, after 7 days. The branching of DNAPL infiltrates further down the aquifer domain. After 15 days, the DNAPL plume reaches up to 9 meter in depth with 2 main branches in the YZ plane. The simulated DNAPL movement along the length of the clay is limited as shown in the XZ section (Figure 6e). This indicates that the accumulation of DNAPL on top of the lower clay is smaller than the required accumulation for DNAPL migration along the length of the clay. Such limitation in accumulation can be attributed to channeling of DNAPL plume along the branches formed along the width of the lower clay. Thus, the model result indicates that the morphology and dimension of the clay plays an important role in determining in the plume migration.

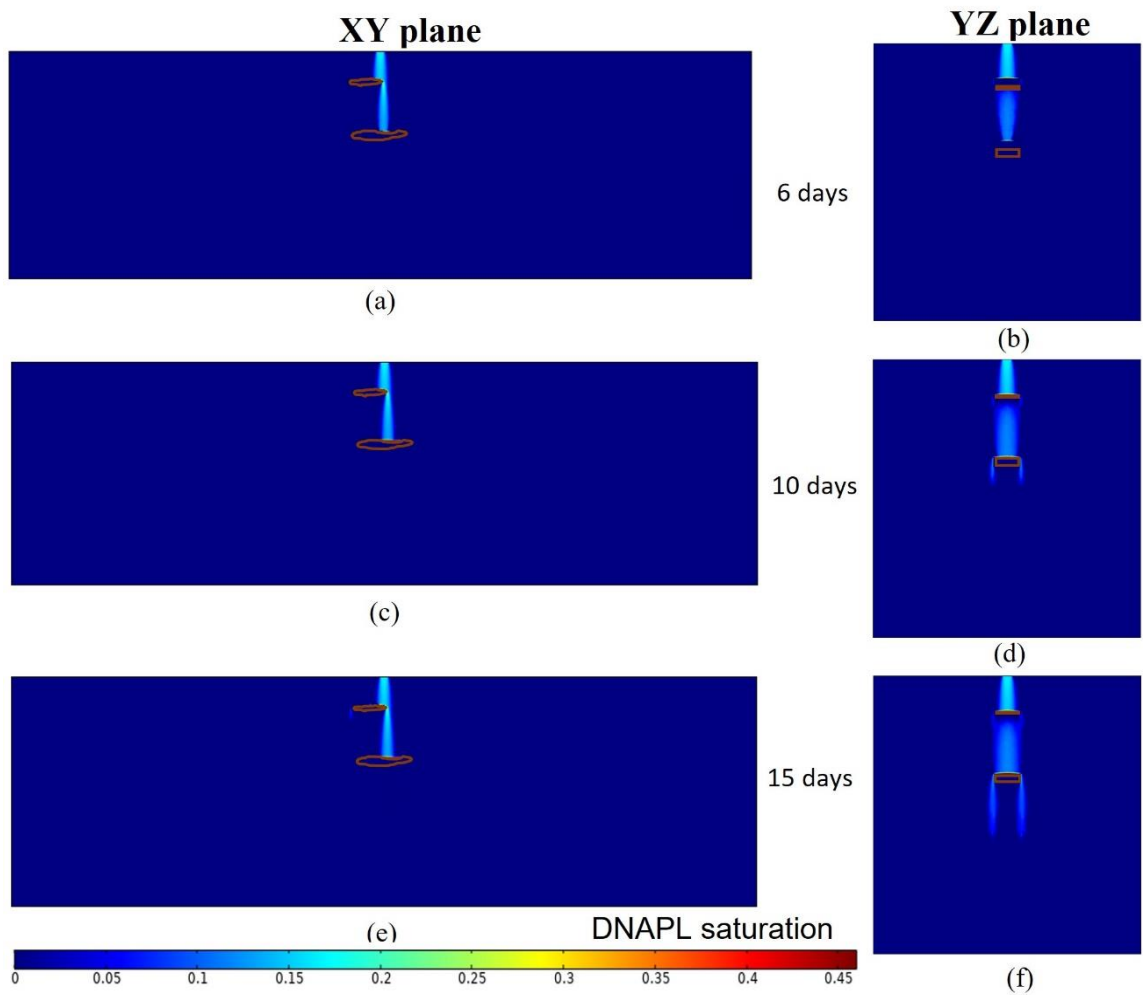

*Figure 6*: Model prediction for DNAPL saturation in the study domain along the XZ plane (a) and YZ plane (b) after 6 days. After 10 days the model predicts further infiltration along the XZ plane (c) and YZ plane (d). After 15 days, when the source condition is removed, the model estimation shows the DNAPL profile along the XZ plane (e) and YZ plane (f).

For a better visualization of DNAPL flow in the aquifer, the iso-surface of DNAPL has been studied. Figure 7a shows the iso-surface with 0.1 DNAPL saturation after 15 days of infiltration in the aquifer. The iso-surface geometry highlights the role of clay in the interruption of vertical movement of DNAPL. The model predicts lateral movement of DNAPL plume on the top of upper clay while the major mass of DNAPL plume follows the original vertical path (Figure 7b). Interestingly, the influence of the lower clay lens on the DNAPL flow path is different from the upper clay lens (Figure 7c). The model predicts a complete halting of the vertical migration of the main DNAPL plume at the top of the lower clay lens. Consequently, two new branches of DNAPL plume are predicted to form along the periphery of the lower clay lens in the model domain. The different flow behavior for the upper and lower clay lenses can be attributed to the different morphology, dimension and position of the clay lenses.

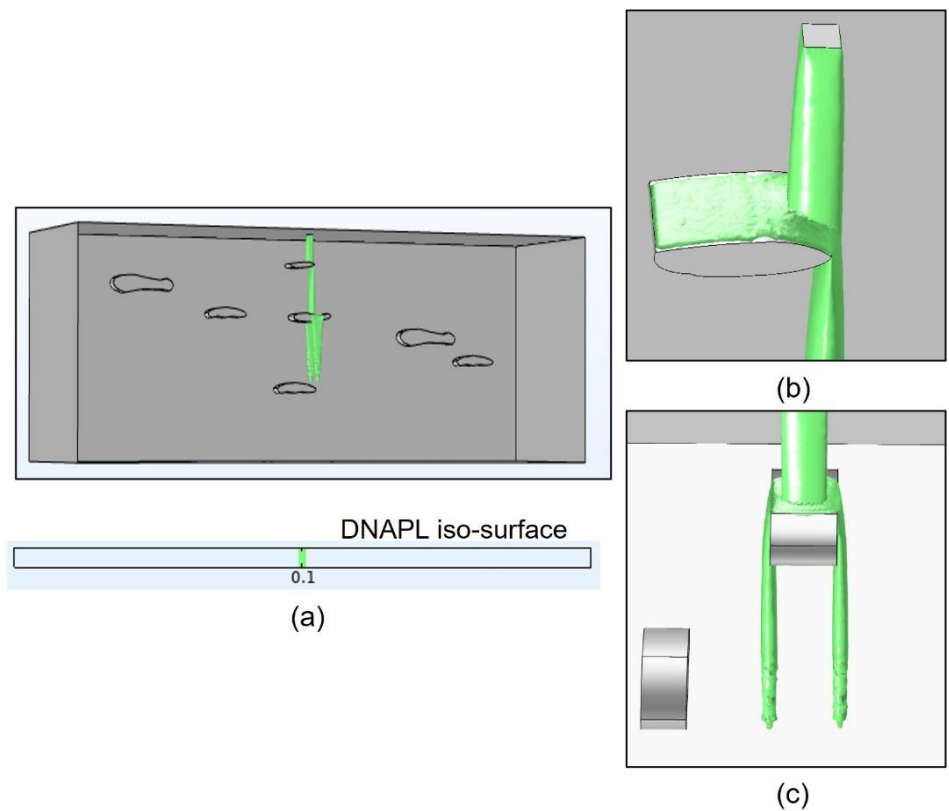

*Figure 7*: Green colored iso-surface of 0.1 DNAPL saturation predicted by model after 15 days of infiltration(a). A closer view for influence of clay lens on the DNAPL movement is shown as an enlarged image along the XZ plane (b) and YZ plane (c).

## 3.2 Prediction of DNAPL migration after source removal

Once the infiltration of the DNAPL from the source has stopped, the model predicts that the DNAPL plume flows down, although with a slower velocity. After 60 days, the flow of DNAPL becomes negligible leading to the formation of stable DNAPL distribution. Figures 8a and 8c show the iso-surface of DNAPL with saturation equal to 0.02 after 1500 hours from different angles of view. The iso-surface represents the model prediction of the formation of a continuous DNAPL ganglia (the interconnected network of residual DNAPL saturation) down to the depth of 5.2 meter below which two vertical branches of DNAPL ganglia occur from the lower clay. The model calculates that the plume front touches the bedrock after 45 days of infiltration leading to the gradual formation of a DNAPL pool on the bedrock. The gradual formation is caused by the slower flow of DNAPL. The DNAPL plume travels 6 meters in 30 days after the source is stopped, making the average velocity 7.5 times slower than the velocity with active infiltration. The decrease of the flow speed can be explained by the drop in the relative permeability due to the decrease of DNAPL saturation. Additionally, as the saturation decreases from the top, the capillary forces start acting upwards due to change in direction of DNAPLs saturation gradient. This leads to the capillary force and gravity force acting in opposite direction. The overall force dynamics results in the formation of a stable DNAPL source zone after 1500 hours. The force dynamics also explains why ganglia have saturation slightly higher than the residual saturation of DNAPL. The final DNAPL source zone has ganglia extending continuously through out the vertical dimension of the simulated aquifer. In the lateral direction, the model predicts the ganglia extending up to 2.5 meter in lateral direction. The negligible movement of DNAPL predicted by model leads to formation of unique plume geometry with ganglia having a DNAPL saturation in the range of 0.02 to 0.05. This agrees with the previous studies (Nambi and Powers, 2003), which reported that the DNAPL ganglia can be in the range of 0.01 to 0.15. This heterogenous distribution of ganglia is such that a low saturation of ganglia is predicted near the water table, while a comparatively higher saturation is predicted with depth, This can be explained by gravity fed motion of finite DNAPL at the source, leading to diminishing saturation starting from the top. The simulated pool formed on the top of clay has decreased from 0.45 to 0.16 after removal of the source. The decrease in DNAPL pool saturation can be attributed to the DNAPL fraction contribution in the vertical motion of the plume. The DNAPL pool formed at the bedrock has the highest saturation of 0.34.

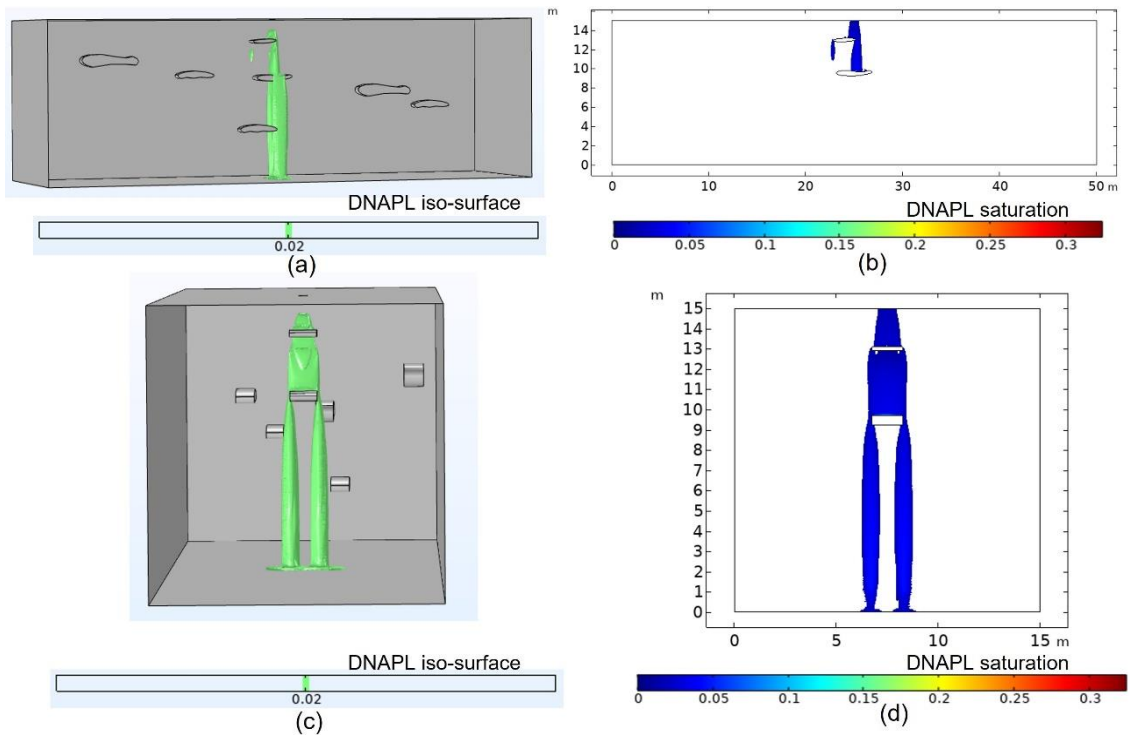

*Figure 8: Model prediction for iso-surface of 0.02 DNAPL saturation in the study domain after 1500 hours viewed along the XZ plane (a) and YZ Plane (c). The range of DNAPL saturation after 1500 hours predicted by model along the reference plane XZ (b) and YZ(d).*

The evolution of contaminant in the simulated aquifer can be studied by plotting the temporal evolution of total mass of DNAPL, ganglia and DNAPL pool (Figure 9). The constant rate of DNAPL infiltration in the simulated domain is reflected as the straight curve shown by orange line for time shorter than 15 days. A total of 500 kg of DNAPL infiltration is simulated by the model in the first fifteen days, and, as the source is removed after 15 days, the total mass of DNAPL remain conserved as represented by flat portion of the orange curve. For the first fifteen days, the simulated DNAPL plume saturation remains higher than 0.16 thus forming a pool of DNAPL along the migration pathway. Therefore, the DNAPL mass represented by green curve shows the predominant fraction of DNAPL mass as DNAPL pool. However, after 15 days, as the source is removed, the gravity leads to reduction of saturation along the migration pathway resulting in formation of ganglia. Then, a rise in the amount of ganglia with a diminishing amount of DNAPL pool is predicted. Finally, as the DNAPL plume reaches the bedrock, a new pool forms slightly increasing its mass proportions. The model estimates a total of 465 kg of ganglia and 32 kg of DNAPL pool at the end of the simulation.

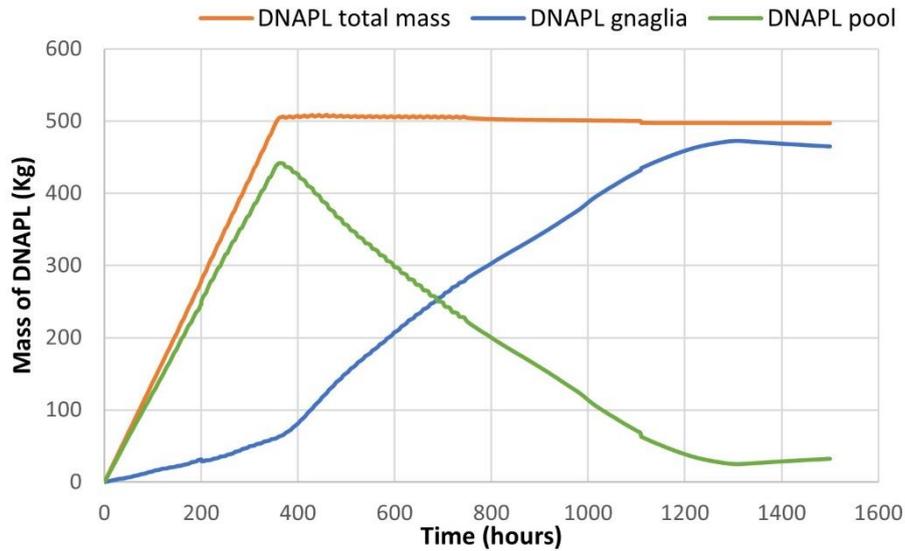

*Figure 9*: Relative masses of DNAPL forms predicted by model after 1500 days. The green curve represents the amount of DNAPL pool, blue curve the amount of ganglia and the orange curve the total amount of DNAPL entered in the simulated domain.

# 4 Implications for the characterisation of polluted sites

The numerical model presented in this study is aimed to predict the flow behavior of DNAPL in a case study. However, its results can be further applied to understand the migration of DNAPL in many other similar aquifers worldwide. Importantly, this kind of predictions are very convenient for the assessment, selection and optimisation of application of remediation technologies.

A first implication from the modelling results is that if the morphology of clay lenses is known, the DNAPL plume geometry can be predicted. This is illustrated in this case study where the model results indicate that the vertical branching depends on both DNAPL plume and the clay size and shape. The preferred pathway for vertical motion of DNAPL below the clay is oriented along the direction with minimum lateral migration of DNAPL over the clay. The lateral migration over the clay becomes negligible as the preferred pathway for vertical migration is developed, resulting in the failure of DNAPL branch formation in other directions. A consequence of this is that the field characterisation can be misleading if not done carefully. During the field study, if a vertically limited DNAPL infiltration in the field is observed, still a large infiltration of DNAPL is possible along the other side of the clay lens. To illustrate this, Figure 10 shows the DNAPL saturation, with a non-uniform formation of branches in the aquifer. The black circle shows a small ganglia formation predicted by model. These ganglia could be the model representation of a small pocket of DNAPL, extending up to 4.2 meter in depth. Therefore, this indicates that tracing the interconnected residual plume requires careful study as the predictive insights from model indicates the possibility of an isolated contamination source in the aquifer.

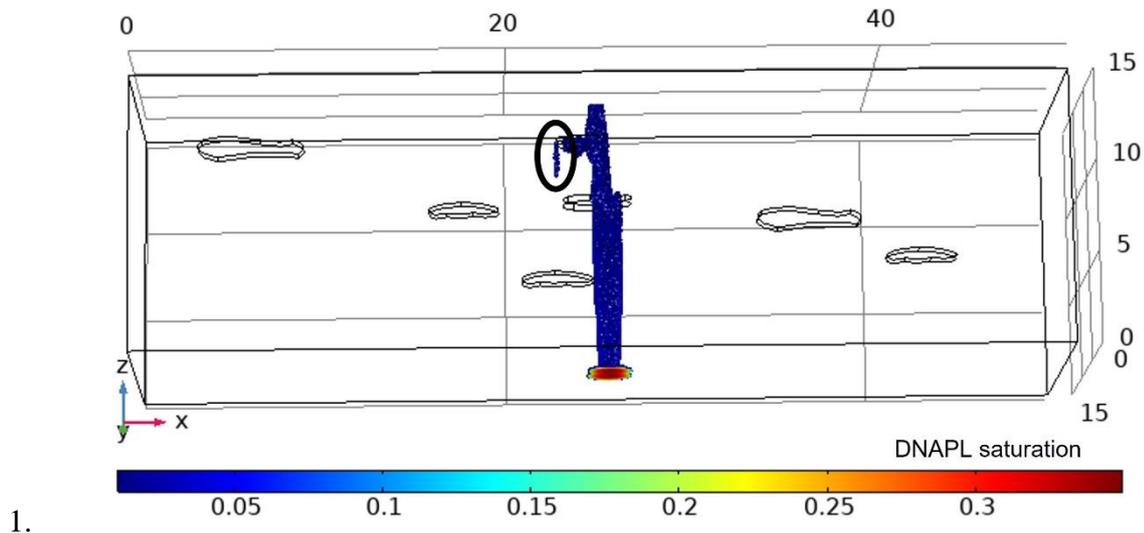

*Figure 10: DNAPL saturation predicted by the model after 1500 hours of DNAPL infiltration in the study domain. The black circle shows the prediction of isolated DNAPL ganglia separated from the main DNAPL plume structure.*

Another interesting output from the numerical model is that the DNAPL accumulated over the clay lens has reduced from 0.45 on day $5^{th}$ to 0.16 after removal of the source. The highest saturation (0.35) of DNAPL occurs at the bedrock. This means that the natural DNAPL attenuation occurs at the shallower level of the aquifer, due to its downward motion leading to accumulation of DNAPL at the bottom. Then, DNAPL is more likely to be present at the bottom part of the aquifer including the top of bedrock in mature fields with a moderate heterogeneity and having history of DNAPL infiltration extended in time. With the maturity of the field, the ganglia dissolve leading its exhaustion (Carey *et al.*, 2014) with only DNAPL pool remaining in the aquifer system. The fate of pool depends on morphology of each impervious clay lenses that the DNAPL plume encounters. For a more flat- or concave-shaped lens surface, the DNAPL pool saturation decreases over time as predicted by the model.

Model results also predicts that the branching of DNAPL plume may or may not re-join. According to the result, DNAPL plume at the shallower clay is infiltrates as a continuous vertical mass. In contrast, the plume branching at the second clay in depth is already discontinuous, and depending on the clay lens geometry, the branches can later re-join. Therefore, at a field scale observation, the DNAPL plume below a clay lens should not be expected to occur only as a discrete form.

The spread of the DNAPL plume predicted by the model is up to 3.7 meter laterally indicating the role of non-permeable layer in DNAPL migration. While the initial DNAPL source inlet has a dimension of 0.5m × 0.5 m width, the lateral movement of DNAPL has led to its widespread distribution. The lateral extension depends on the number of impermeable layers it encounters. The extent of vertical migration of DNAPL, however, depends on the heterogeneity of the aquifer as well as on the mass of infiltrated DNAPL. For a small mass infiltrating down into the aquifer with higher heterogeneity, the DNAPL plume may not reach the bedrock of the aquifer. Such high heterogeneity of aquifer can lead to accumulation of DNAPL on the less permeable layer, thus leading to the limited vertical movement of DNAPL. Additionally, the higher heterogeneity in an aquifer can lead to higher lateral migration of DNAPL in the aquifer system.

# 5 Conclusions

The site remediation strategy in DNAPL polluted areas requires a lot of effort and resources of characterisation, owing to complexity and inaccessibility of the contaminant distribution. The success of the remediation strategy is often limited by lack of detailed understanding of the contaminant and hydrogeological properties of the impacted soils and aquifers. While, the information of contamination-aquifer dynamics can be studied at the well location, a continuous and consistent contamination behavior on both temporal and spatial scale is often missing. In this article, a field scale two-phase flow model is developed for a better understanding of the migration behaviour of the DNAPL in a hypothetical aquifer with meso-scale heterogeneities. The model results based on the fluid and hydrogeological features of a hypothetical aquifer demonstrates that the complexity of DNAPL migration is highly governed by the heterogeneities in the porous media. While the DNAPL plume motion is predicted to be mainly vertical in the homogeneous sand, the presence of clay lenses causes lateral migration of DNAPL. The predicted flow behavior results in a complex DNAPL source zone architecture involving formation of several DNAPL branches. The simulation successfully quantifies the varying amount of DNAPL pool and ganglia formed over the time. The correct numerical implementation of Interface Condition is verified by studying the predicted capillary pressure and saturation at the sand-clay interfaces. The numerical integration of fluid mass conservation equations along with interface condition enables the simulation to illustrate the role of low permeable clays within the aquifer in DNAPL flow dynamics and immobilization. Overall, the numerical tool demonstrates that the predictive insight of DNAPL behavior can help understand the extent of contamination and develop suitable remediation strategy.